\documentclass{article}
\usepackage{amsmath, amssymb, graphics}
\usepackage[pdflatex]{graphicx}
\usepackage{amsthm}
\usepackage{amscd}
\usepackage{caption}

\newtheorem{thm}{Theorem}[section]
\newtheorem{df}{Definition}[section]

\newtheorem{prp}{Proposition}[section]

\newtheorem{cor}{Corollary}[section]
  
\newcommand{\mathsym}[1]{{}}
\newcommand{\unicode}[1]{{}}

\makeatletter
    
    \@addtoreset{equation}{section}
  \makeatother
  
\begin{document}
\title{Mathematical Analysis of Melodies: Slope and Discrete Fr\'{e}chet distance}
\author{Fumio HAZAMA\\Tokyo Denki University\\Hatoyama, Hiki-Gun, Saitama JAPAN\\
e-mail address:hazama@mail.dendai.ac.jp\\Phone number: (81)49-296-2911}
\date{\today}
\maketitle
\thispagestyle{empty}

\begin{abstract}
A directed graph, called an {\it M-graph}, is attached to every melody. Our chief concern in this paper is to investigate (1) how the positivity of the slope of the M-graph is related to {\it singability} of the melody, (2) when the M-graph has a symmetry, and (3) how we can detect a similarity between two melodies. For the third theme, we introduce the notion of {\it transposed discrete Fr\'{e}chet distance}, and show its relevance in the study of similarity detection among an arbitrary set of melodies.
\end{abstract}

\section{Introduction}
In the article [2], the authors introduced a method of attaching a graph to an arbitrary melody. We call it here the {\it M-graph} of the melody. By using the M-graph of melody as a main ingredient, we investigate in this paper (1) how the positivity of the slope of the M-graph is related to {\it singability} of the melody, (2) when the M-graph has a symmetry, and (3) how we can detect a similarity between two melodies. Accordingly we divide the paper into three parts. In the first part we focus on the slope of the M-graph, which is defined by the  method of least squares, and investigate how the slope is related to  musical characteristics of the original melody. For example, among melodies which are composed of six notes {C4, D4, E4, F4, G4, A4} and begin with C4, the largest slope is attained by (C4, D4, E4, F4, G4, A4) with slope 0.986 and the smallest slope is attained by (C4, A4, D4, G4, E4, F4)) with slope -0.729. One can see that the latter is harder to sing than the former. Through the analysis of several data including this example, we will show that the positivity of the slope is strongly related to its singability. One word of caution: We do not assert that positivity of the slope is related to its goodness. For example, the first phrase of the most famous nocturne (in E$\flat$ major) by Chopin has a (slightly) negative slope -0.089, but cannot be claimed that it is a bad melody accordingly. We see, however, that all of the fifteen other nocturnes by Chopin have positive slopes (see Table 8). We also consider how the slope of the M-graph is changed under transposition, inversion, and retrograde of the original melody. In the second part of the paper, we investigate how a symmetry of the M-graph is reflected to the character of the melody. We invite the reader to have a look at Fig. 2, which is the M-graph of the basic twelve-tone row of the string quartet Op. 28 by Webern. This amazing example leads us to the main theorem (Theorem 2.1) of the second part, which characterizes the melodies with symmetric M-graph in terms of a certain arithmetic property. In the third part of the paper we propose a distance, called {\it transposed discrete Fr\'{e}chet distance}, and show its relevance for similarity detection through several examples. The data in the final subsection come from the author's questionnaire to the students in a class on discrete geometry. The national anthem of Israel, "Twinkle, twinkle, little star", and the Japanese classical song "Kojo no Tsuki", which constitute the nearest cluster, are found to be sung simultaneously and quite harmoniously. This surprise motivated him to write this paper.\\

\section{Slope of M-graph}
\subsection{Definition of M-graph}
 In order to express a melody by a definite sequence of integers, we let C4 (middle C) correspond to 0, C$\#$4 to 1, and so on. In this way we can associate a sequence of integers with each melody. For example the melody "C4, D4, F4, E4", which is the main theme of the fourth movement of the Jupiter symphony by Mozart, corresponds to the sequence "0, 2, 5, 4". From now on we identify a melody of finite length with the sequence of integers of finite length which is constructed by this rule. Furthermore,to any sequence $\mathbf{a}=(a_1,a_2,\cdots,a_n)$ of integers, we attach a sequence of points $\mathbf{p}=(p_1,p_2,\cdots,p_{n-1})$ with $p_i\in\mathbf{R}^2\hspace{1mm}(1\leq i\leq n-1)$ by the following rule:
 \begin{eqnarray*}
 p_1=(a_1,a_2), p_2=(a_2,a_3),\cdots, p_{n-1}=(a_{n-1},a_n).
 \end{eqnarray*}
Let $G(\mathbf{a})=(V(\mathbf{a}),E(\mathbf{a}))$ be the directed graph with the set of vertices
\begin{eqnarray*}
V(\mathbf{a})=(p_1,p_2,\cdots,p_{n-1}),
\end{eqnarray*}
and the set of edges
\begin{eqnarray*}
E(\mathbf{a})=\{(p_1,p_2),(p_2,p_3), \cdots, (p_{n-2},p_{n-1})\}.
\end{eqnarray*}
We call $G(\mathbf{a})$ the {\it M-graph} associated to the melody $\mathbf{a}$. ("M" stands for \underline{m}elody.) When $\mathbf{a}=(0,2,5,4)$, for example, its M-graph $G(\mathbf{a})$ is depicted as follows:

\newpage
\begin{figure}[ht]
\centering
\includegraphics{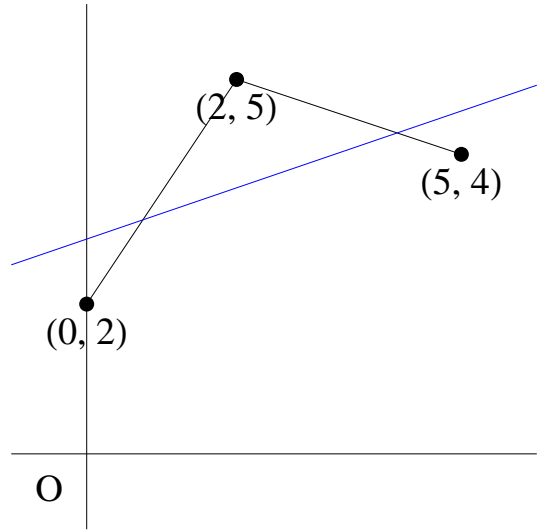}
  \caption{M-graph of Jupiter}
 \end{figure}

\noindent
The line which cuts through the M-graph in this figure is obtained by the least squares fitting. Its slope will be referred as {\it the slope} of the melody, and denoted by $s(\mathbf{a})$. In this case we see that $s(\mathbf{a})=0.342$.\\

\noindent
Remark. A formula for the slope in the method of least squares will be recalled in Proposition 1.1.

\subsection{Distribution of slopes of M-graphs}
We will show that there exists a correlation between the slope of a melody and its {\it singability}. Let us look at the set 
\begin{eqnarray*}
M_4=\{(0,2,4,5), (0,2,5,4), (0,4,2,5), (0,4,5,2), (0,5,2,4), (0,5,4,2)\},
\end{eqnarray*}
which collects all the melodies consisting of C4, D4, E4, F4 beginning with C4. The slopes of these are computed as follows:

\newpage
\begin{table}[ht]
 \begin{center}
  \begin{tabular}{llc}
   \hline
   melody & note names & slope  \\
   \hline \hline
   (0,2,4,5) & (C,D,E,F) & 0.750\\
   (0,2,5,4) & (C,D,F,E) & 0.342\\
   (0,4,2,5) & (C,E,D,F) & -0.500\\
   (0,4,5,2) & (C,E,F,D) & -0.214\\
   (0,5,2,4) & (C,F,D,E) & -0.605\\
   (0,5,4,2) & (C,F,E,D) & -0.357\\
   \hline
  \end{tabular}
  \end{center}
  \caption{Slopes of four-tone melodies}
\end{table}

\noindent
Notice here that our friend (C,D,F,E) has the second highest slope among the melodies in $M_4$, and that the other melodies, except the simplest melody (C,D,E,F), have negative slopes. In order to understand what is going on, we take next the set $M_5$ of melodies consisting of C4, D4, E4, F4, G4 beginning with C4. The top three melodies with largest slope and the bottom three with smallest slope are tabulated below:
\begin{table}[ht]
 \begin{center}
  \begin{tabular}{lllc}
   \hline
   ranking & melody & note names & slope  \\
   \hline \hline
   1st & (0,2,4,5,7) & (C,D,E,F,G) & 0.915\\
   2nd & (0,2,5,4,7) & (C,D,F,E,G) & 0.576\\
   3rd & (0,2,4,7,5) & (C,D,E,G,F) & 0.467\\
   \hline
  \end{tabular}
  \end{center}
  \caption{Largest three slopes}
\end{table}\\
\begin{table}[ht]
 \begin{center}
  \begin{tabular}{lllc}
   \hline
   ranking & melody & note names & slope  \\
   \hline \hline
   -1st & (0,7,2,5,4) & (C,G,D,F,E) & -0.655\\
   -2nd & (0,7,2,4,5) & (C,G,D,E,F) & -0.617\\
   -3rd & (0,7,4,5,2) & (C,G,E,F,D) & -0.538\\
   \hline
  \end{tabular}
  \end{center}
  \caption{Smallest three slopes}
\end{table}

\noindent
The next table shows the top three and the worst three of slopes among melodies which consists of six notes C4, D4, E4, F4, G4, A4 and begins with C4:

\newpage
\begin{table}[ht]
 \begin{center}
  \begin{tabular}{lllc}
   \hline
   ranking & melody & note names & slope  \\
   \hline \hline
   1st & (0,2,4,5,7,9) & (C,D,E,F,G,A) & 0.98630\\
   2nd & (0,2,5,4,7,9) & (C,D,F,E,G,A) & 0.81507\\
   3rd & (0,2,4,7,5,9) & (C,D,E,G,F,A) & 0.64384\\
   3rd & (0,4,2,5,7,9) & (C,E,D,F,G,A) & 0.64384\\
   \hline
  \end{tabular}
  \end{center}
  \caption{Largest three slopes}
\end{table}
\begin{table}[ht]
 \begin{center}
  \begin{tabular}{lllc}
   \hline
   ranking & melody & note names & slope  \\
   \hline \hline
   -1st & (0,9,2,7,4,5) & (C,A,D,G,E,F) & -0.72932\\
   -2nd & (0,7,4,5,2,9) & (C,G,E,F,D,A) & -0.72603\\
   -3rd & (0,9,2,7,5,4) & (C,A,D,G,F,E) & -0.69925\\
   \hline
  \end{tabular}
  \end{center}
  \caption{Smallest three slopes}
\end{table}

\noindent
As the reader may notice in these examples, melodies with large (positive) slope tend to be easy to sing and those with small (negative) slope are hard to sing. The table below describes the numbers of melodies with positive, negative, or zero slope in each category:

\begin{table}[ht]
 \begin{center}
  \begin{tabular}{lccc}
   \hline
  constituent & positive & negative & zero  \\
   \hline \hline
  \{C,D,E,F,G\} & 8 & 16 & 0\\
  \{C,D,E,F,G,A\} & 45 & 75 & 0\\
  \{C,D,E,F,G,A,B\} & 262 & 457 & 1\\
   \hline
  \end{tabular}
  \end{center}
  \caption{Distribution of slopes}
\end{table}

\noindent
We notice that, in each category, the number of melodies with negative slope is about twice the number of those with positive slope. Therefore we may assert that composers choose instinctively melodies with positive slope, which constitue rather a minor part in the world of melodies, in order to make their works singable ones. \\

Keeping these observations in mind, we examine the slopes of actual melodies composed by two great composers, Schumann  and Chopin. Table 7 shows the slopes of the first phrases of the sixteen songs in "Dichterliebe" by Schumann:

\newpage
\begin{table}[ht]
 \begin{center}
  \begin{tabular}{lllllllll}
   \hline
  No. & 1 & 2 & 3 & 4 & 5 & 6 & 7 & 8  \\
   \hline
  slope & 0.183 & 0.302 & 0.951 & 0.553 & 0.545 & 0.712 & 0.438 & 0.584 \\
   \hline
  \end{tabular}
  \end{center}
\end{table}
\begin{table}[ht]
 \begin{center}
  \begin{tabular}{lllllllll}
   \hline
  No. & 9 & 10 & 11 & 12 & 13 & 14 & 15 & 16  \\
   \hline
  slope & 0.691 & 0.543 & 0.656 & 0.316 & 0.929 & -0.572 & 0.450 &0.666 \\
   \hline
  \end{tabular}
  \end{center}
  \caption{Distribution of slopes in the song cycle Dichterliebe}
\end{table}

\noindent
Among these songs, only one song has a negative slope. It is the fourteenth song, titled  "Alln\"{a}chtlich im Traume seh' ich dich," whose slope is -0.572. The fact that, when we listen to the song cycle as a whole, we feel a certain soothing effect at this 14-th song, might be related to the negativity of its slope. \\

The following table shows the slopes of all the Nocturnes composed by Chopin:\\

\newpage
\begin{table}[ht]
 \begin{center}
  \begin{tabular}{lllllllll}
   \hline
  No. & 1 & 2 & 3 & 4 & 5 & 6 & 7 & 8  \\
   \hline
  slope & 0.980 & -0.089 & 0.371 & 0.508 & 0.860 & 0.667 & 0.496 & 0.677 \\
   \hline
  \end{tabular}
  \end{center}
\end{table}
\begin{table}[ht]
\begin{center}
  \begin{tabular}{lllllllll}
   \hline
  No. & 9 & 10 & 11 & 12 & 13 & 14 & 15 & 16  \\
   \hline
  slope & 0.419 & 0.641 & 0.673 & 0.650 & 0.970 & 0.520 & 0.293 &0.099 \\
   \hline
  \end{tabular}
  \end{center}
  \caption{Distribution of slopes in Nocturnes}
\end{table}

\noindent
Among these nocturnes, only the second one has a negative slope. This nocturne also has a kind of soothing effect, which might be one of the reasons why this is widely regarded as the most popular Nocturne by Chopin. On the other hand the largest and the second largest slopes are attained by the first one (in B$\flat$ minor) and 13-th one (in C minor), respectively. Both pieces move us (or at least the author) with their distinctive deep sorrow.\\

\subsection{Slopes under transformations}
In this subsection we consider what occurs to the slope of a melody if it is transposed, inverted, or reversed.\\

For an arbitrary melody $\mathbf{x}=(x_1,\cdots,x_{n+1})$ and for any $t\in\mathbb{Z}$, let $\mathbf{x}+t=(x_1+t,\cdots,x_{n+1}+t)$, the transposition by $t$. Furthermore we denote the inversion $(-x_1,\cdots,-x_{n+1})$ by $\mathbf{x}^i$, and the retrograde $(x_{n+1},\cdots,x_1)$ by $\mathbf{x}^r$. Here we recall the formula for the slope of a point data based on the method of least squares:

\begin{prp}
Let $P$ denote a set of points $(x_1,y_1),\cdots,(x_N,y_N)$ on $\mathbb{R}^2$. Then the slope $s(P)$ obtained through the method of least squares is given by the formula
\begin{eqnarray}
s(P)=\frac{N\sum_{i=1}^Nx_iy_i-\sum_{i=1}^Nx_i\sum_{i=1}^Ny_i}{N\sum_{i=1}^Nx_i^2-\left(\sum_{i=1}^Nx_i\right)^2}.
\end{eqnarray}
\end{prp}

\noindent
(I) Transposition: The M-graph of transposition $M(\mathbf{x}+t)$ consists of the points $(x_i+t, x_{i+1}+t)\hspace{1mm}(1\leq i \leq n)$. Hence we have $M(\mathbf{x}+t)=M(\mathbf{x})+(t,t)$, namely all of the points in $M(\mathbf{x+t})$ are translations of the ones in $M(\mathbf{x})$ by one and the same point $(t,t)$. Therefore by the very definition of the method of least squares we have the following:

\begin{prp}
For any melody $\mathbf{x}$, we have
\begin{eqnarray*}
s(M(\mathbf{x}+t))=s(M(\mathbf{x})).
\end{eqnarray*}
\end{prp}

\noindent
Remark. One can prove this by a direct computation of the slope on the right hand side by using the formula (1.1).\\

\noindent
(II) Inversion: Since the numerator and the denominator of the right hand side of (1.1) are homogeneous polynomials of degree two in the variables $x_i, y_i\hspace{1mm}(1\leq i\leq N)$, both of them are invariant under the transformation $x_i\mapsto -x_i, y_i\mapsto -y_i\hspace{1mm}(1\leq i\leq N)$. Hence we have the following:

\begin{prp}
For any melody $\mathbf{x}$, we have
\begin{eqnarray*}
s(M(\mathbf{x}^i))=s(M(\mathbf{x})).
\end{eqnarray*}
\end{prp}

\noindent
Combination of Proposition 1.2 and 1.3 yields the following:

\begin{cor}
For any melody $\mathbf{x}$ and for any $t\in\mathbb{Z}$, we have
\begin{eqnarray*}
s(M(\mathbf{x}^i+t))=s(M(\mathbf{x})).
\end{eqnarray*}
\end{cor}

\noindent
For example, if $\mathbf{x}=$(A4, C5, B4, A4, E5)$=(9,12,11,9,16)$ (Paganini), then $\mathbf{x}^i+17=(8,5,6,8,1)=$(A$\flat$4, F4, G$\flat$4, A$\flat$4, D$\flat$4) (Rachmaninov). It follows from Corollary 1.1 that their slopes coincide. Actually one can see that $s(M(\mathbf{x}))=-1.333<0$, and passingly that $s(M($B$\flat$3, C4, D$\flat$4, A$\flat$3))$=-1.071<0$, but that  their concatenation satisfies $s(M($A$\flat$4, F4, G$\flat$4, A$\flat$4, D$\flat$4, B$\flat$3, C4, D$\flat$4, A$\flat$3))=0.668 . Thus Rachmaninov composed this fascinating melody with positive slope by combining the two parts with negative slope.\\

\noindent
(III) Retrograde: It turns out to be essential to deal with the numerator and the denominator of the slope separately. Accordingly we set for any melody $\mathbf{x}=(x_1,\cdots,x_{n+1}),$
\begin{eqnarray*}
N(\mathbf{x})&=&n\sum_{i=1}^nx_ix_{i+1}-\sum_{i=1}^nx_i\sum_{i=2}^{n+1}x_i,\\
D(\mathbf{x})&=&n\sum_{i=1}^nx_i^2-\left(\sum_{i=1}^nx_i\right)^2,
\end{eqnarray*}
which are obtained by setting $y_i=x_{i+1}\hspace{1mm}(i=1,\cdots,n)$ and $N=n$ in (1.1). Let us put $\mathbf{x}^r=(x_1',\cdots,x_{n+1}')$ so that $x_i'=x_{n+2-i}\hspace{1mm}(1\leq i\leq n+1)$. First we look at the numerator $N(\mathbf{x})$:

\begin{prp}
For any melody $\mathbf{x}$, we have
\begin{eqnarray*}
N(\mathbf{x}^r)=N(\mathbf{x}).
\end{eqnarray*}
\end{prp}

\noindent
{\it Proof}. This can be proved by the following straightforward computation:
\begin{eqnarray*}
N(\mathbf{x}^r)&=&n\sum_{i=1}^nx'_ix'_{i+1}-\sum_{i=1}^nx'_i\sum_{i=2}^{n+1}x'_i\\
&=&n\sum_{i=1}^nx_{n+2-i}x_{(n+2)-(i+1)}-\sum_{i=1}^nx_{n+2-i}\sum_{i=2}^{n+1}x_{n+2-i}\\
&&(by\hspace{1mm}letting\hspace{1mm}i'=n+1-i)\\ 
&=&n\sum_{i'=1}^nx_{i'+1}x_{i'}-\sum_{i'=1}^nx_{i'+1}\sum_{i'=0}^{n-1}x_{i'+1}\\
&=&N(\mathbf{x}).
\end{eqnarray*}
\qed\\

\noindent
The denominator is, however, not invariant under the retrograde transformation:

\begin{prp}
For any melody $\mathbf{x}$, we have
\begin{eqnarray*}
D(\mathbf{x}^r)=D(\mathbf{x})
\end{eqnarray*}
if and only if
\begin{eqnarray*}
x_1=x_{n+1}\hspace{1mm}\mbox{or}\hspace{1mm}(n+1)(x_{n+1}+x_1)=2\sum_{i=1}^{n+1}x_i.
\end{eqnarray*}
\end{prp}

\noindent
{\it Proof}. We compute the difference $D(\mathbf{x}^r)-D(\mathbf{x})$:
\begin{eqnarray*}
D(\mathbf{x}^r)-D(\mathbf{x})&=&\left(n\sum_{i=1}^n(x'_i)^2-\left(\sum_{i=1}^nx'_i\right)^2\right)\\
&&-\left(n\sum_{i=1}^nx_i^2-\left(\sum_{i=1}^nx_i\right)^2\right)\\
&=&\left(n\sum_{i=2}^{n+1}x_i^2-\left(\sum_{i=2}^{n+1}x_i\right)^2\right)\\
&&-\left(n\sum_{i=1}^nx_i^2-\left(\sum_{i=1}^nx_i\right)^2\right)\\
&=&n(x_{n+1}^2-x_1^2)+\left(\left(\sum_{i=1}^nx_i\right)^2-\left(\sum_{i=2}^{n+1}x_i\right)^2\right)\\
&=&n(x_{n+1}-x_1)(x_{n+1}+x_1)\\
&&+\left(\sum_{i=1}^nx_i-\sum_{i=2}^{n+1}x_i\right)\left(\sum_{i=1}^nx_i+\sum_{i=2}^{n+1}x_i\right)\\
&=&n(x_{n+1}-x_1)(x_{n+1}+x_1)\\
&&+(x_1-x_{n+1})\left(2\sum_{i=1}^{n+1}x_i-x_1-x_{n+1}\right)\\
&=&(x_{n+1}-x_1)\left((n+1)(x_{n+1}+x_1)-2\sum_{i=1}^{n+1}x_i\right).
\end{eqnarray*}
Hence the assertion follows.
\qed\\

\noindent
Combining Propositions 1.4 and 1.5, we have the following:

\begin{cor}
When a melody $\mathbf{x}$ begins and ends with one and the same note, $\mathbf{x}$ and its retrograde have the same slope.
\end{cor}

\subsection{Locality of the slope function}
In this subsection we will see that the slope of a melody is determined by the slopes of its parts. More precisely we show the following:

\begin{prp}
For any melody $\mathbf{x}=(x_1,\cdots,x_{n+1})$ and for any $k$ with $1\leq k\leq n-1$, let $\mathbf{x}^k$ denote the triple $(x_k,x_{k+1},x_{k+2})$ and let $s_k=s(M(\mathbf{x}^k))$. Then the slope $s(M(\mathbf{x}))$ of the whole melody is a rational function in $s_1,\cdots,s_{n-1}$.
\end{prp}

\noindent
{\it Proof}. It follows from Corollary 1.1 that $s(M(\mathbf{x}))=s(M(\mathbf{x}-x_1)$. Hence if we put $y_i=x_{i+1}-x_i\hspace{1mm}(1\leq i\leq n)$, then we have
\begin{eqnarray*}
s(M(\mathbf{x}))&=&s(M(0,y_1,y_1+y_2,\cdots,\sum_{i=1}^ny_i)\\
&=&\frac{N(0,y_1,y_1+y_2,\cdots,\sum_{i=1}^ny_i)}{D(0,y_1,y_1+y_2,\cdots,\sum_{i=1}^ny_i)},
\end{eqnarray*}
the rightmost side being the ratio of homogeneous quadratic polynomials in $y_1,\cdots,y_n$. Hence dividing the numerator and the denominator by $y_1^2$, we see that $s(M(\mathbf{x}))$ is a quotient of quadratic polynomials in $y_2/y_1, y_3/y_1, \cdots, y_n/y_1$. Hence it is a rational function in $y_2/y_1, y_3/y_2, \cdots, y_n/y_{n-1}$. On the other hand, we see that 
\begin{eqnarray*}
s_k=s(M(\mathbf{x}^k))&=&s(\{(x_k,x_{k+1}), (x_{k+1},x_{k+2})\})\\
&=&\frac{x_{k+2}-x_{k+1}}{x_{k+1}-x_k}\\
&=&\frac{y_{k+1}}{y_k}
\end{eqnarray*}
holds for any $k$ with $1\leq k\leq n-1$. This completes the proof.\qed\\

\noindent
Example 3.1. When $\mathbf{x}=(x_1,x_2,x_3,x_4)$ is a melody of length 4, we have
\begin{eqnarray*}
s(M(\mathbf{x}))&=&\frac{3(x_1x_2+x_2x_3+x_3x_4)-(x_1+x_2+x_3)(x_2+x_3+x_4)}{3(x_1^2+x_2^2+x_3^2)-(x_1+x_2+x_3)^2}\\
&=&\frac{y_2^2+2y_1y_2+y_1y_3+2y_2y_3}{2(y_1^2+y_1y_2+y_2^2)}\\
&=&\frac{(y_2/y_1)^2+2(y_2/y_1)+(y_3/y_1)+2(y_2/y_1)(y_3/y_1)}{2(1+(y_2/y_1)+(y_2/y_1)^2)}\\
&=&\frac{s_1^2+2s_1+s_1s_2+2s_1^2s_2}{2(1+s_1+s_1^2)}.
\end{eqnarray*}

\noindent
Remark. For an arbitrary finite set of points $P$ in the plane, the slope $s(P)$ is \underline{not} necessarily a function in the slopes of consecutive segments. For example, let $P=\{(0,0), (1,0), (2,1)\}$. Then the slopes of consecutive segments are 0 and 1, and the whole slope is computed to be 
\begin{eqnarray*}
s(P)=\frac{3\cdot2-3\cdot 1}{3\cdot 5-3^2}=\frac{1}{2}.
\end{eqnarray*}
On the other hand if we put $P'=\{(0,0), (1,0), (3,2)\}$, then the consecutive slopes are 0 and 1, and hence the local slopes coincide with those of $P$. The whole slope, however, turns out to be
\begin{eqnarray*}
s(P')=\frac{3\cdot6-4\cdot 2}{3\cdot 10-4^2}=\frac{5}{7}.
\end{eqnarray*}
Thus the slope $s(P)$ is not generally a function of local slopes.\\

\noindent
\section{Symmetry of M-graphs}
In this section we investigate for what kind of melodies their associated M-graphs have reflective symmetries.\\

For a line $\ell$ in the plane $\mathbb{R}^2$, let $ref_{\ell}:\mathbb{R}^2\rightarrow\mathbb{R}^2$ denote the reflection with the line $\ell$ as a set of fixed points. We introduce the following:

\begin{df}
For any melody $\mathbf{x}=(x_1,\cdots,x_{n+1})$, let $M(\mathbf{x})=(p_1,\cdots,p_n)$ be its M-graph so that $p_i=(x_i,x_{i+1})$ for $i=1,\cdots,n$. The melody $\mathbf{x}$ is said to have a reflective symmetry if there exists a line $\ell$ such that $ref_{\ell}(p_i)=p_{n+1-i}$ holds for any $i\in[1,n]$.
\end{df}

\noindent
For example, when $\mathbf{x}=(0,1,2,\cdots,n)$, then one can see that $\mathbf{x}$ has a reflective symmetry with respect to the line $y=-x+n$. We want to characterize the set of melodies with reflective symmetry. For this purpose we need a transformation formula. When $\ell$ is defined by the equation $y=ax+b$, we have
\begin{eqnarray}
ref_{\ell}(x,y)=\left(\frac{(1-a^2)x+2ay-2ab}{1+a^2},\frac{2ax-(1-a^2)y+2b}{1+a^2}\right).
\end{eqnarray}
On the other hand, when $\ell$ is parallel to the $y$-axis and hence defined by the equation $x=c$, we have
\begin{eqnarray*}
ref_{\ell}(x,y)=(-x+2c,y).
\end{eqnarray*}

\noindent
In the present paper we restrict our attention to the melodies without repetition, namely those with pairwise distinct entries.
First we deal with the melodies of even length $2n$, and we denote a general melody by indexing it as
\begin{eqnarray*}
\mathbf{x}=(x_{-n},x_{-(n-1)},\cdots,x_{-1},x_1,\cdots,x_{n-1},x_n).
\end{eqnarray*}
 This will ease our description of an inductive argument. We start with the case $n=2$.\\

\begin{prp}
A melody $\mathbf{x}=(x_{-2},x_{-1},x_1,x_2)$ without repetition has a reflective symmetry if and only if the following condition is satisfied:\\
\begin{eqnarray}
&&{\rm (I)}\hspace{2mm}x_2=-x_{-2}+x_{-1}+x_1,\\
&&\mbox{or}\nonumber\\
&&{\rm (II)}\hspace{2mm}x_2=x_{-2}-x_{-1}+x_1.
\end{eqnarray}
The respective axis of symmetry is given by
\begin{eqnarray}
&&{\rm (I)}\hspace{2mm}y=-x+x_{-1}+x_1,\\
&&{\rm (II)}\hspace{2mm}y=\frac{x_{-2}-x_1}{x_{-2}-2x_{-1}+x_1}x-\frac{(x_{-1}-x_1)(x_{-2}+x_1)}{x_{-2}-2x_{-1}+x_1}.
\end{eqnarray}
\end{prp}

\noindent
{\it Proof}. Let $\ell$ be the axis of symmetry. Then the following two conditions must be met:
\begin{eqnarray}
ref_{\ell}(x_{-2},x_{-1})&=&(x_1,x_2),\\
ref_{\ell}(x_{-1},x_1)&=&(x_{-1},x_1). 
\end{eqnarray}

\noindent
By our assumption we have $x_{-1}\neq x_2$, and hence the condition (2.6) implies that $\ell$ is not parallel to the $y$-axis. Let $y=ax+b$ be its defining equation. It follows from the formula (2.1) that the condition (2.6) leads us to the following simultaneous equation
\begin{eqnarray}
  \left\{
    \begin{array}{l}
      \frac{(1-a^2)x_{-2}+2ax_{-1}-2ab}{1+a^2} = x_1 \\
      \frac{2ax_{-2}-(1-a^2)x_{-1}+2b}{1+a^2} = x_2
    \end{array}
  \right.
\end{eqnarray}
By multiplying $1+a^2$ on both sides of these equations, we have
\begin{eqnarray}
  \left\{
    \begin{array}{l}
     (x_{-2}+x_1)a^2-2x_{-1}a+2ab-x_{-2}+x_1=0 \\
    (x_{-1}-x_2)a^2+2x_{-2}a+2b-x_{-1}-x_2=0
    \end{array}
  \right.
\end{eqnarray}
By subtracting the first equation from the second equation multiplied by $a$, we obtain
\begin{eqnarray*}
(x_{-1}-x_2)a^3+(x_{-2}-x_1)a^2+(x_{-1}-x_2)a+(x_{-2}-x_1)=0,
\end{eqnarray*}
namely we have
\begin{eqnarray*}
(a^2+1)((x_{-1}-x_2)a+(x_{-2}-x_1))=0.
\end{eqnarray*}
Since $a$ is a real number and $x_{-1}-x_2\neq 0$, we see that
\begin{eqnarray}
a=\frac{x_{-2}-x_1}{-x_{-1}+x_2}.
\end{eqnarray}
Inserting this expression into the second equation of (2.9), we find that
\begin{eqnarray}
b=\frac{x_{-2}^2+x_{-1}^2-x_1^2-x_2^2}{2(x_{-1}-x_2)}.
\end{eqnarray}
Furthermore the condition (2.7) with these values for $a$ and $b$ is expressed as the equalities
\begin{eqnarray*}
&&x_{-2}^3-x_{-1} x_{-2}^2-x_1 x_{-2}^2+x_{-1}^2 x_{-2}-x_1^2 x_{-2}-x_2^2 x_{-2}+2 x_1 x_2
   x_{-2}\\
&&+x_{-1}^3+x_1^3+x_{-1} x_1^2+x_{-1} x_2^2+x_1 x_2^2-x_{-1}^2 x_1-2 x_{-1}^2 x_2-2 x_1^2
   x_2\\
&&=(x_{-2}^2-2 x_1 x_{-2}+x_{-1}^2+x_1^2+x_2^2-2 x_{-1} x_2)x_{-1},\\
&&x_{-1}^3-2 x_{-2} x_{-1}^2+x_1
   x_{-1}^2-x_2 x_{-1}^2+x_{-2}^2 x_{-1}-x_1^2 x_{-1}-x_2^2 x_{-1}\\
&&+2 x_{-2} x_2 x_{-1}+x_1^3+x_2^3-2 x_{-2}
   x_1^2-x_1 x_2^2+x_{-2}^2 x_1-x_{-2}^2 x_2+x_1^2 x_2\\
&&=(x_{-2}^2-2 x_1
   x_{-2}+x_{-1}^2+x_1^2+x_2^2-2 x_{-1} x_2)x_1.
\end{eqnarray*}
These equations are factored, somewhat miraculously, as
\begin{eqnarray*}
(x_{-2}-x_1)(x_{-2}-x_{-1}-x_1+x_2)(x_{-2}-x_{-1}+x_1-x_2)&=&0,\\
(x_{-1}-x_2)(x_{-2}-x_{-1}-x_1+x_2)(x_{-2}-x_{-1}+x_1-x_2)&=&0.
\end{eqnarray*}
Since $x_{-1}-x_2\neq 0$, the second equation implies that  $x_{-2}-x_{-1}-x_1+x_2=0$ or $x_{-2}-x_{-1}+x_1-x_2=0$, and both alternatives satisfy the first equation. Hence we have
\begin{eqnarray}
{\rm (I) }\hspace{1em}x_2=-x_{-2}+x_{-1}+x_1,\mbox{ {\rm or (II)} }x_2=x_{-2}-x_{-1}+x_1.
\end{eqnarray}
In case of (I), the slope $a$ and the $y$-intercept $b$ are found through (2.10) and (2.11) to be
\begin{eqnarray*}
a&=&\frac{x_{-2}-x_1}{-x_{-1}+(-x_{-2}+x_{-1}+x_1)}=\frac{x_{-2}-x_1}{-x_{-2}+x_1}=-1,\\
b&=&\frac{x_{-2}^2+x_{-1}^2-x_1^2-(-x_{-2}+x_{-1}+x_1)^2}{2(x_{-1}-(-x_{-2}+x_{-1}+x_1))}\\
&=&\frac{-2x_1^2+2(x_{-2}x_{-1}+x_{-2}x_1-x_{-1}x_1)}{2(x_{-2}-x_1)}\\
&=&\frac{2(x_{-2}-x_1)(x_{-1}+x_1)}{2(x_{-2}-x_1)}\\
&=&x_{-1}+x_1.
\end{eqnarray*}
This shows that the axis of symmetry in this case is given by (2.4), and the reflection map is given by
\begin{eqnarray*}
ref_{\ell}:(x,y)\mapsto (-y+x_{-1}+x_1, -x+x_{-1}+x_1).
\end{eqnarray*}
Therefore the condition (2.2) is also sufficient for the reflective symmetry of $\mathbf{x}$. In case of (II), a similar computation based on (2.10) and (2.11) shows that (2.5) holds true. Furthermore we notice the following interesting phenomenon in this case: the triangle $p_1p_2p_3$ is a isoceles right triangle with $\angle p_1p_2p_3=90^{\circ}$. For we have
\begin{eqnarray*}
p_2-p_1&=&(x_{-1}-x_{-2},x_1-x_{-1}),\\
p_3-p_2&=&(x_1-x_{-1},x_2-x_1)=(x_1-x_{-1},x_{-2}-x_{-1}),
\end{eqnarray*}
which are transversal and have equal lengths. Therefore the melody $\mathbf{x}=(x_{-2},x_{-1},x_1,x_{-2}-x_{-1}+x_1)$ has a reflective symmetry with the bisector of $\angle p_1p_2p_3$ as the axis of symmetry. 
It follows that the condition (2.3) is also sufficient for the melody $\mathbf{x}=(x_{-2},x_{-1},x_1,x_2)$ to have a reflective symmetry. This completes the proof. \qed\\

\noindent
The following corollary can be deduced easily from Proposition 2.1, but it will facilitate our inductive argument later:

\begin{cor}
If a melody $\mathbf{x}=(x_{-2},x_{-1},x_1,x_2)$ has a reflective symmetry, then we have
\begin{eqnarray*}
x_2-x_1=x_{-1}-x_{-2}\mbox{ or }x_2-x_1=x_{-2}-x_{-1}.
\end{eqnarray*}
\end{cor}

Next we consider the melodies with six notes.
\begin{prp}
A melody $\mathbf{x}=(x_{-3},x_{-2},x_{-1},x_1,x_2,x_3)$ without repetition has a reflective symmetry if and only if 
\begin{eqnarray*}
x_{-i}+x_i\mbox{ is constant for }i=1,2,3.
\end{eqnarray*}
The axis of symmetry is given by
\begin{eqnarray*}
y=-x+x_{-1}+x_1,
\end{eqnarray*}
and the reflection map  is given by
\begin{eqnarray*}
(x,y)\mapsto (-y+x_{-1}+x_1,-x+x_{-1}+x_1).
\end{eqnarray*}
\end{prp}

\noindent
{\it Proof}. Since the submelody $(x_{-2},x_{-1},x_1,x_2)$ also has a reflective symmetry, we are in the cases (I) or (II) in Proposition 2.1. \\

\noindent
(I) The case when $x_2=-x_{-2}+x_{-1}+x_1$ : The reflection in this case is given by
\begin{eqnarray*}
ref_{\ell}(x,y)=(-y+x_{-1}+x_1,-x+x_{-1}+x_1),
\end{eqnarray*}
and hence we must have
\begin{eqnarray*}
ref_{\ell}(x_{-3},x_{-2})=(-x_{-2}+x_{-1}+x_1,-x_{-3}+x_{-1}+x_1)=(x_2,x_3)
\end{eqnarray*}
Therefore we have $x_{-3}+x_3=x_{-1}+x_1$.\\

\noindent
(II) The case when $x_2=x_{-2}-x_{-1}+x_1$: Let $q_{24}$ (resp. $q_{15}$) denote the midpoint of $p_2p_4$ (resp. $p_1p_5$). Recalling that
\begin{eqnarray*}
&&p_1=(x_{-3},x_{-2}),\hspace{2mm}p_2=(x_{-2},x_{-1}),\\
&&p_4=(x_1,x_{-2}-x_{-1}+x_1),\hspace{2mm}p_5=(x_{-2}-x_{-1}+x_1,x_3),
\end{eqnarray*}
we have
\begin{eqnarray*}
&&q_{24}=(\frac{x_{-2}+x_1}{2},\frac{x_{-2}+x_1}{2}),\\
&&q_{15}=(\frac{x_{-3}+x_{-2}-x_{-1}+x_1}{2},\frac{x_{-2}+x_3}{2}).
\end{eqnarray*}
Hence we have
\begin{eqnarray}
&&\overrightarrow{p_3q_{24}}=(\frac{x_{-2}-2x_{-1}+x_1}{2},\frac{x_{-2}-x_1}{2}),\\
&&\overrightarrow{q_{24}q_{15}}=(\frac{x_{-3}-x_{-1}}{2},\frac{-x_1+x_3}{2}).
\end{eqnarray}
Note that these vectors are nonzero by our assumption. Furthermore, since these two vectors have the same direction with the axis of symmetry, the equality
\begin{eqnarray}
\overrightarrow{p_3q_{24}}=k\overrightarrow{q_{24}q_{15}}
\end{eqnarray}
holds for some $k\in\mathbb{R}^*$. Since the submelody $(x_{-3},x_{-2},x_2,x_3)$ must have the same axis of symmetry as the one for $(x_{-2},x_{-1},x_1,x_2)$, it follows from Corollary 2.1 that we necessarily have
\begin{eqnarray}
x_3-x_2=x_{-3}-x_{-2},
\end{eqnarray}
or
\begin{eqnarray}
x_3-x_2=x_{-2}-x_{-3}.
\end{eqnarray}
Accordingly we divide our argument further into two cases.\\

\noindent
Case (II.A): The case when $x_3-x_2=x_{-3}-x_{-2}$. Since we are in the case when $x_2=x_{-2}-x_{-1}+x_1$, we have
\begin{eqnarray*}
x_3-(x_{-2}-x_{-1}+x_1)=x_{-3}-x_{-2},
\end{eqnarray*}
which implies that
\begin{eqnarray*}
x_3-x_1=x_{-3}-x_{-1}.
\end{eqnarray*}
Hence we have
\begin{eqnarray*}
\overrightarrow{q_{24}q_{15}}=(\frac{x_{-3}-x_{-1}}{2},\frac{-x_1+x_3}{2})=\frac{x_{-3}-x_{-1}}{2}(1,1).
\end{eqnarray*}
This implies by (2.15) and (2.14) that the $x$-coordinate and the $y$-coordinate of $\overrightarrow{p_3q_{24}}$ must coincide, and hence it follows from (2.13) that
\begin{eqnarray*}
x_{-2}-2x_{-1}+x_1=x_{-2}-x_1.
\end{eqnarray*}
This implies that $x_{-1}=x_1$, which contradicts to our assumption. \\

\noindent
Case (II.B): The case when $x_3-x_2=x_{-2}-x_{-3}$: Note that the vector $\overrightarrow{p_1p_5}$ and the vector $\overrightarrow{p_3q_{24}}$ are transversal by symmetry assumption. Since in our case
\begin{eqnarray*}
\overrightarrow{p_1p_5}&=&(x_2,x_3)-(x_{-3},x_{-2})=(x_2-x_{-3})\cdot (1,1)\neq (0,0),\\
\end{eqnarray*}
it follows from (2.13) that the inner product of $\overrightarrow{p_1p_5}$ and $\overrightarrow{p_3q_{24}}$, which is equal to
\begin{eqnarray*}
&&\frac{x_2-x_{-3}}{2}((x_{-2}-2x_{-1}+x_1)+(x_{-2}-x_1))\\
&=&(x_2-x_{-3})(x_{-2}-x_{-1}),
\end{eqnarray*}
must become zero. This, however, contradicts to our assumption. Hence both of Case (II.A) and Case (II.B) cannot occur, and the proof is completed. \qed\\

Now we can generalize Proposition 2.2 to an arbitrary melody of even length:

\begin{thm}
For any integer $n\geq 3$, a melody $\mathbf{x}=(x_{-n},\cdots,x_{-1},x_1,\cdots,x_n)$ without repetition of length $2n$ has a reflective symmetry if and only if it satisfies the following condition:
\begin{eqnarray*}
{\rm (I)}_n:\hspace{1mm}x_{-i}+x_i\mbox{ is constant for }i=1,\cdots,n.
\end{eqnarray*}
When this condition is met, the axis $\ell$ of symmetry is the line defined by $y=-x+x_{-1}+x_1$, and the reflection map is given by
\begin{eqnarray*}
ref_{\ell}: (x,y)\mapsto (-y+x_{-1}+x_1,-x+x_{-1}+x_1).
\end{eqnarray*}
\end{thm}

\noindent
{\it Proof}. We prove this by induction on $n$. When $n=3$, this is Proposition 2.2 itself. When $n\geq 4$, suppose that  a melody $\mathbf{x}=(x_{-n},\cdots,x_{-1},x_1,\cdots,x_n)$ has a reflective symmetry. Then the submelody $(x_{-(n-1)},\cdots,x_{-1},x_1,\cdots,x_{n-1})$ also has a reflective symmetry. Then by the induction hypothesis, the assertion ${\rm (I)}_{n-1}$ holds true. Then we have
\begin{eqnarray}
x_{-(n-1)}+x_{n-1}=x_{-1}+x_1.
\end{eqnarray}
Since the reflection in this case is given by
\begin{eqnarray*}
ref_{\ell}(x,y)=(-y+x_{-1}+x_1,-x+x_{-1}+x_1),
\end{eqnarray*}
we must have
\begin{eqnarray*}
ref_{\ell}(x_{-n},x_{-(n-1)})&=&(-x_{-(n-1)}+x_{-1}+x_1,-x_{-n}+x_{-1}+x_1)\\
&=&(x_{n-1},x_n).
\end{eqnarray*}
The equality of the first entries is assured by (2.18), and that for the second entries is equivalent to
\begin{eqnarray}
x_{-n}+x_n=x_{-1}+x_1.
\end{eqnarray}
Hence the assertion ${\rm (I)}_n$ holds. Conversely, suppose that the condition ${\rm (I)}_n$ holds, and let $\ell$ be the line defined by $y=-x+x_{-1}+x_1$. Then the reflection map with the axis of symmetry $\ell$ is given by
\begin{eqnarray*}
ref_{\ell}: (x,y)\mapsto (-y+x_{-1}+x_1,-x+x_{-1}+x_1).
\end{eqnarray*}
If follows that, for any $k\in [1,n]$, we have
\begin{eqnarray*}
ref_{\ell}(x_{-k},x_{-(k-1)})&=&(-x_{-(k-1)}+x_{-1}+x_1,-x_{-k}+x_{-1}+x_1)\\
&=&(x_{k-1},x_k),
\end{eqnarray*}
which shows that the melody $\mathbf{x}$ has a reflective symmetry with the axis of symmetry $\ell$. This completes the proof. \qed\\

\noindent
Remark. For a melody of odd length without repetition, we can show the following result: When $n\geq 2$, a melody $\mathbf{x}=(x_{-n},\cdots,x_{-1},x_0,x_1,\cdots,x_n)$ without repetition has a reflective symmetry if and only if $x_{-i}+x_i=2x_0$ for any $i\in [1,n]$. This can be proved in a similar way to that for Theorem 2.1, so we omit the proof.\\

\noindent
Example. 2.1. Webern based his string quartet Op. 28 on the following row:
\begin{eqnarray*}
M_W=(7,6,9,8,12,13,10,11,15,14,17,16).
\end{eqnarray*}
Amazingly, this melody of length twelve turns out to satisfy the condition ${\rm (I)}_{6}$ in Theorem 2.1. Hence it must have a reflective symmetry. We illustrate below the M-graph of its transposition
\begin{eqnarray*}
M_{W'}=(1,0,3,2,6,7,4,5,9,8,11,10).
\end{eqnarray*}
(Note that transposing does not change the reflective property of the original melody.)

\begin{figure}[ht]
\centering
\includegraphics{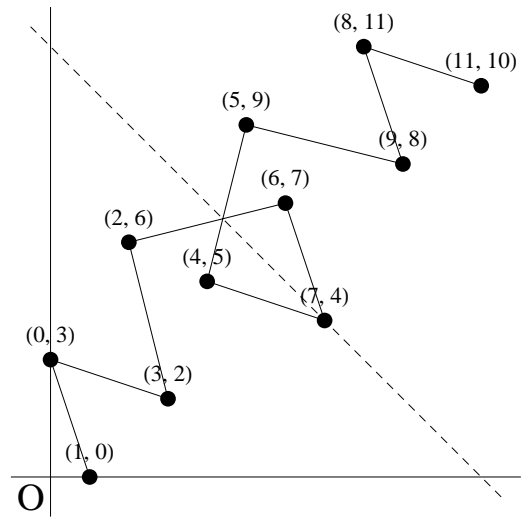}
  \caption{String quartet Op. 28 by Webern}
 \end{figure}

\noindent
The dashed line is the axis of symmetry of $M_{W'}$ and is defined by the equation $y=-x+11$.\\

\noindent
Example. 2.2. The following row is used in Ode to Napoleon Op.41 by Sh\"{o}nberg:
\begin{eqnarray*}
M_S=(1,0,4,5,9,8,3,2,6,7,11,10).
\end{eqnarray*}
Here again we are surprised that this melody satisfies the condition ${\rm (I)}_{6}$ in Theorem 2.1. Hence it has a reflective symmetry:

\newpage
\begin{figure}[ht]
\centering
\includegraphics{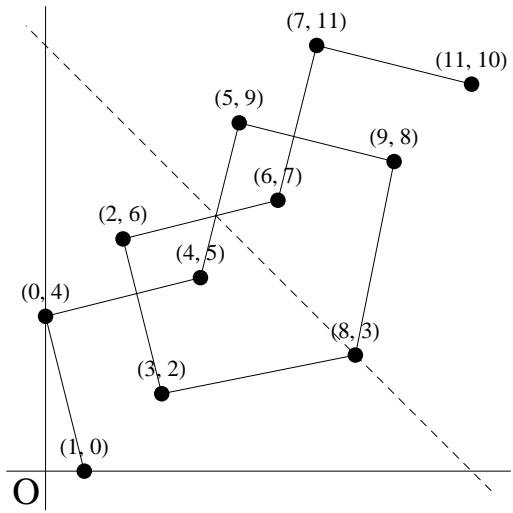}
  \caption{Ode to Napoleon Op.41 by Sh\"{o}nberg}
 \end{figure}

\noindent
The dashed line is the axis of symmetry of $M_S$ and is defined by the equation $y=-x+11$.\\

\noindent
It is needless to say that every twelve-tone row does not have a reflective symmetry. Thus these two composers arrived at the above symmetrical rows through their musical intellect and instinct.

\section{Transposed Discrete Fr\'{e}chet distance}
In this section we introduce the notion of {\it transposed discrete Fr\'{e}chet distance}, abbreviated as TDFD. This is based on the discrete Fr\'{e}chet distance, abbreviated as DFD. We will show the relevance of TDFD for similarity detection among a given set of melodies. \\

\subsection{Definition of DFD and TDFD}
First we recall the definition of DFD for the convenience of the reader. (See [1], [3] for details.) Let $P=(p_1,p_2,\cdots,p_n)$ and $Q=(q_1,q_2,\cdots,q_m)$ be a pair of sequences of points in $\mathbb{R}^2$. A {\it coupling} $L$ between $P$ and $Q$ is a sequence
\begin{eqnarray*}
(p_{a_1},q_{b_1}), (p_{a_2},q_{b_2}), \cdots, (p_{a_k},q_{b_k})
\end{eqnarray*}

\noindent
of distinct pairs from $P\times Q$ such that $a_1=b_1=1, a_k=n, b_k=m$, and for any $i=1,\cdots,k-1$, we have $a_{i+1}=a_i$ or $a_{i+1}=a_i+1$, and $b_{i+1}=b_i$ or $b_{i+1}=b_i+1$. The {\it length} $||L||$ of the coupling $L$ is defined by
\begin{eqnarray*}
||L||={\rm max}_{i=1,\cdots,k}d(p_{a_i},q_{b_i}),
\end{eqnarray*}

\noindent
where $d(*,*)$ denotes the Euclid distance on $\mathbf{R}^2$. The discrete Fr\'echet distance $d_F(P,Q)$ between the sequences of points $P$ and $Q$ is defined to be
\begin{eqnarray*}
d_{F}(P,Q)={\rm min}\{||L||;\text {$L$ is a coupling between $P$ and $Q$}\}
\end{eqnarray*}

\noindent
Intuitively this can be defined as follows. A man is walking a dog on a leash. The man can move on the points in the sequence $P$, and the dog in the sequence $Q$, but backtracking is not allowed. The discrete Fr\'echet distance $d_F(P,Q)$ is the length of the shortest leash that is sufficient for traversing both sequences. For a pair of melodies $\mathbf{a}, \mathbf{b}$, we define the discrete Fr\'echet distance $d_F(\mathbf{a},\mathbf{b})$ to be $d_F(V(\mathbf{a}),V(\mathbf{b}))$. Furthermore taking into account the fact that any transposition of a melody does not change its essential feature, we define the transposed discrete Fr\'echet distance $d_{F}^{tr}(\mathbf{a},\mathbf{b})$ by the following rule:
\begin{eqnarray}
d_{F}^{tr}(\mathbf{a},\mathbf{b})=\min_{t\in\mathbb{Z}}d_F(\mathbf{a},\mathbf{b}+t)
\end{eqnarray}

\noindent
where $\mathbf{b}+t$ denotes the transposed melody $(b_1+t,\cdots,b_m+t)$. In an actual computation of $d_{F}^{tr}(\mathbf{a},\mathbf{b})$, we can choose a bound $B$ such that the minimum on the right hand side of (2.1) lies in $[-B,B]$. \\

\noindent
Example 3.1. Let $\mathbf{a}_1=(0,2,4,5,2,2,0)$ and $\mathbf{b}_1=(0,2,5,2,1)$. The point sequences which correspond to these melodies are 
\begin{eqnarray*}
V(\mathbf{a}_1)=(p_1, p_2, \cdots, p_6),
\end{eqnarray*}
with
\begin{eqnarray*}
p_1=(0,2), p_2=(2,4), p_3=(4,5), p_4=(5,2), p_5=(2,2), p_6=(2,0)),
\end{eqnarray*}
and
\begin{eqnarray*}
V(\mathbf{b}_1)=(q_1, q_2, \cdots, q_4),
\end{eqnarray*}
with
\begin{eqnarray*}
q_1=(0,2), q_2=(2,5), q_3=(5,2), q_4=(2,1).
\end{eqnarray*}

\noindent
The coupling $L$ which attains the minimum of $||L||$ is found to be
\begin{eqnarray*}
L=((p_1, q_1), (p_2, q_2), (p_3, q_2), (p_4, q_3), (p_5,q_4),(p_6,q_4))
\end{eqnarray*}
with $||L||=2$. Note that when the man Paul (for $P$) goes to $p_3$, his dog Queen (for $Q$) must remain at $q_2$, because if Queen moves to $q_3$, then $d(p_3,q_3)=\sqrt{10}>2=d(p_3,q_2)$. I recommend the reader to take a walk with Queen several times, then he will be convinced that the above coupling $L$ is the best choice.\\

\noindent
Example 3.2. Let 
\begin{eqnarray*}
\mathbf{a}_2&=&(0,2,4,5,7)=(C4, D4, E4, F4, G4),\\
\mathbf{b}_2&=&(2,9,7,6,4)=(D4, A4, G4, F\#4, E4).
\end{eqnarray*}
 The discrete Fr\'{e}chet distances between $\mathbf{a}_2$ and $\mathbf{b}_2+t$ with $t=-5, -4, \cdots, 0, 1$ are tabulated as follows:\\

\begin{table}[ht]
 \begin{center}
  \begin{tabular}{cl}
   \hline
  $t$ & $d_F(\mathbf{a}_2, \mathbf{b}_2+t)$\\
   \hline \hline
 -5 & 8.944\\
 -4 & 7.616\\
 -3 & 6.325\\
 -2 & 5.099\\
 -1 & 6.083\\
 0 & 7.280\\
 1 & 8.544\\
    \hline
  \end{tabular}
  \end{center}
  \caption{discrete Fr\'{e}chet distances}
\end{table}

\noindent
It follows that $d_{F}^{tr}(\mathbf{a}_2,\mathbf{b}_2)=d_F(\mathbf{a}_2,\mathbf{b}_2-2)=5.099$. This seems to be natural, since the melody $\mathbf{a}_2$ is in C major, the melody $\mathbf{b}_2$ in D major, and $C4-D4=0-2=-2$. The next example, however, shows us that the situation is not so simple.\\

\noindent
Example 3.3. Let 
\begin{eqnarray*}
\mathbf{a}_3&=&(0,2,4,5,7)=(C4, D4, E4, F4, G4),\\
\mathbf{b}_3&=&(0,4,7,12)=(C4, E4, G4, C5).
\end{eqnarray*}
The discrete Fr\'{e}chet distances between $\mathbf{a}_3$ and $\mathbf{b}_3+t$ with $t=-5, -4, \cdots, 0, 1$ are tabulated as follows:\\

\newpage
\begin{table}[ht]
\begin{center}
  \begin{tabular}{cl}
   \hline
  $t$ & $d_F(\mathbf{a}_3, \mathbf{b}_3+t)$\\
   \hline \hline
 -5 & 5.831\\
 -4 & 4.472\\
 -3 & 3.162\\
 -2 & 3.000\\
 -1 & 4.123\\
 0 & 5.385\\
 1 & 6.708\\
    \hline
  \end{tabular}
  \end{center}
  \caption{discrete Fr\'{e}chet distances}
\end{table}

\noindent
It follows that $d_{F}^{tr}(\mathbf{a}_3,\mathbf{b}_3)=d_F(\mathbf{a}_3,\mathbf{b}_3-2)=3.000$. This time both melodies are in C major, but they require a transposition by -2. Indeed the arithmetic mean of the entries in the melody $\mathbf{a}_3$ is 3.600, that for $\mathbf{b}_3$ is 5.750, and their difference is equal to $-2.15\approx -2$. \\

\noindent
In Example 3.2, the arithmetic mean of $\mathbf{a}_2$ is 3.600, that of $\mathbf{b}_2$ is 5.600, and their difference is equal to -2. This together with Example 3.3 shows that the relevance of difference of the arithmetic means of two melodies when we compute the transposed Fr\'{e}chet distance.@These phenomena lead us to consider the DTFD's between a melody and its permutations. Note that in this case their arithmetic means are one and the same. \\

\noindent
Example 3.4. Let us fix $\mathbf{a}_4=(0,2,4)=(C2, D4, E4)$ and let $\mathbf{b}_4$ runs in the set of permutations of $\{0,2,4\}$. The following table displays the values of $t$ for which $d_F(\mathbf{a}_4,\mathbf{b}_4+t)$ attains the minimum:\\

\newpage
\begin{table}[ht]
 \begin{center}
  \begin{tabular}{llcc}
   \hline
  $\mathbf{a}_4$ & $\mathbf{b}_4$ & $t$ with minimum $d_F(\mathbf{a}_4, \mathbf{b}_4+t)$ & distance\\
   \hline \hline
 (0,2,4) & (0,2,4) & 0 & 0\\
 (0,2,4) & (0,4,2) & 0 & 2.828\\
 (0,2,4) & (2,0,4) & 0 & 2.828\\
 (0,2,4) & (2,4,0) & 1 & 4.243\\
 (0,2,4) & (4,0,2) & -1 & 4.243\\
 (0,2,4) & (4,2,0) & 0 & 4.000\\
    \hline
  \end{tabular}
  \end{center}
  \caption{transposed discrete Fr\'{e}chet distances}
\end{table}

\noindent
These examples teach us a lesson that the difference of the arithmetic means is a tentative value for us to find what value of $t$ gives us the minimum of TFD.

\subsection{Cluster analysis based on TDFD}
In this subsection we analyze the cluster structure of some instances of melodies by using TDFD.\\

As samples we choose several national anthems. The following table shows the names of countries, their national anthems in terms of numbers, and their slopes:

\begin{table}[ht]
 \begin{center}
  \begin{tabular}{lllc}
   \hline
   & name of country & national anthem & slope\\
   \hline \hline
 1 & Austria & (12,10,9,10,12,14,12,12,10,10) & 0.460\\
 2 & Bulgaria & (4,9,9,11,12,11,9,4,9,9,11,12,11,9) & 0.257\\
 3 & Canada & (7,10,10,3,5,7,9,10,12,5) & 0.110\\
 4 & China & (7,11,14,14,16,14,11,7,14,14,14,11,7) & 0.285\\
 5 & Germany & (7,9,11,9,12,11,9,6,7,16,14,12,11,9,11,7,14) & 0.131\\
 6 & Hungary & (2,3,5,10,5,3,2,7,5,3,2,0,2,3) & 0.359\\
 7 & Israel & (0,2,3,5,7,7,7,8,7,8,12,7,5,5,5,3,5,3,2,0,2,3,0) & 0.743\\
 8 & Japan & (2,0,2,4,7,4,2,4,7,9,7,9,14,11,9,7) & 0.729\\
 9 & Morocco & (10,12,10,7,8,10,3,5,7,7,8,0,7,8,5) & 0.165\\
 10 & New Zealand & (7,6,7,2,11,11,9,7,4,12,2,11,9,7,6,4,2) & -0.197\\
    \hline
  \end{tabular}
  \end{center}
  \caption{slopes of national anthems}
\end{table}
 
\noindent
For these melodies we apply clustering by using group average method. As a result, we find that the melody 1 (Austria) and the melody 6 (Hungary) is the closest pair. Furthermore the following table reveals a fascinating fact:

\newpage
\begin{table}[ht]
 \begin{center}
  \begin{tabular}{cl}
   \hline
  $t$ & $d_F("Austria", "Hungary"+t)$\\
   \hline \hline
 4 & 7.211\\
 5 & 5.831\\
 6 & 4.472\\
 7 & 3.606\\
 8 & 4.123\\
 9 & 5.385\\
 10 & 6.708\\
    \hline
  \end{tabular}
  \end{center}
  \caption{DFD between "Austria " and transposed "Hungary"}
\end{table}

\noindent
The arithmetic mean of "Austria" is equal to 11.100 and that of "Hungary" is 3.714, and hence the difference is $7.386\approx 7$, which coincides the value of $t$ giving the minimum of DFD under transposition. Actually, "Austria" is in F major, and "Hungary" is in B$\flat$ major. Hence our TDFD detects the difference of their arithmetic means as well as the difference F$4-$B$\flat3=5-(-2)=7$, which is required to transpose "Hungary" to "Austria". Moreover, though they are in different time, one is in three-four time, the other in four-four, both melodies can be played at the same time quite harmoniously. \\

Next we consider how cluster structure changes if we add some other melodies to the above 10 national anthems. Let us choose "Twinkle Twinkle Little Star" as the 11-th melody:
\begin{eqnarray*}
11: {\rm Twinkle}=(0, 0, 7, 7, 9, 9, 7, 5, 5, 4, 4, 2, 2, 0): {\rm slope} = 0.690
\end{eqnarray*}

Here appears a new nearest pair (7: Israel, 11:Twinkle) with $d_F^{tr}(7,11)=3.606$, which is equal to TDFD between 1 and 6. Amazingly enough, the melodies 7 and 11 can be sung harmoniously under the condition that Twinkle is transposed to C minor. \\

Futhermore we add to the samples a Japanese song called "Kojo no Tsuki (Moon over the Ruined Castle)" as the 12-th melody:
\begin{eqnarray*}
12: {\rm Kojo}=(6, 6, 11, 13, 14, 13, 11, 7, 7, 6, 5, 6): {\rm slope} = 0.762
\end{eqnarray*}

Here appears a cluster (7: Israel, 11:Twinkle, 12:Kojo) where $d_F^{tr}(7,12)=d_F^{tr}(11,12)=2.828$. We are surprised again to find that these three melodies can be sung harmoniously if Twinkle and Kojo are transposed to C minor. \\

\subsection{Motivating example}
This subsection explains how the author came across the cluster (7: Israel, 11:Twinkle, 12:Kojo). In the fall term in 2013 he gave lessons in the discrete Fr\'{e}chet distance at his university, and in a class he sent out to the students questionnaire about their most favorite musics. However almost all of the 25 answers which they supplied were Japanese song, the above three melodies happened to be contained in them. By clustering of 25 melodies, he came across to a cluster of the three melodies as well as another cluster of "The Moldau" and "Ievan Polkka". This pair can be sung simultaneously too!.@These surprises in the result of questionaire motivated the author to study in this article the usefulness of the transposed discrete Fr\'echet distance.\\
\newpage
\noindent
References\\
$[1]$ T. Eiter, H. Mannila. Computing discrete Fr\'{e}chet distance. Technical Report CD-TR 94/64, Information Systems Department, Technical University of Vienna, 1994.\\
$[2]$ G. G\"{u}nd\"{u}z, U. G\"{u}nd\"{u}z, The mathematicas analysis of the structure of some songs, Physica A, $\mathbf{357}$(2005), 565-592.\\
$[3]$ A. Mosig, M. Clausen, Approximately matching polygonal curves with respect to the Fr\'{e}chet distance, Comput. Geom. $\mathbf{30}$(2005) 113-127.

 \end{document}